\def\cite{} 
\def\edful{1}
\def\forful{2}
\def\ggt{3}
\def\gpv{4}
\def\gkp{5}
\def\sgb{6}
\def\lsv{7}
\def\mirsk{8}
\def\papst{9}
\def\schrij{10}
\def\spenc{11}

\magnification\magstep1
\baselineskip12pt
\parskip3pt plus 1pt

\def\Sigmait{{\mit\Sigma}}
\def\Sbar{{\overline S}}
\def\xbar{{\bar{x}}}
\def\Sigbar{{\overline\Sigmait}}
\def\ldt{\mathrel{.\,.}}
\def\rar{\rightarrow}
\def\pfbox
  {\hbox{\hskip 3pt\lower2pt\vbox{\hrule
  \hbox to 7pt{\vrule height 7pt\hfill\vrule}
  \hrule}}\hskip3pt}
\def\proof{\noindent{\bf Proof. \ }}
\setbox0=\hbox{[10]\enspace} \newdimen\bibhang \bibhang=\wd0
\def\bib{\par\noindent\hangindent\bibhang}
\def\rbar{\relbar\!\!\relbar}

\centerline{\bf Two-way Rounding}
\smallskip
\centerline{\sl Donald E. Knuth}
\centerline{\sl Computer Science Department}
\centerline{\sl Stanford University}
\centerline{\sl Stanford, CA 94305}

\bigskip
{\narrower\smallskip\noindent
{\bf Abstract.} 
Given $n$ real numbers $0\leq x_1,\ldots,x_n<1$ and a permutation~$\sigma$ of
$\{1,\ldots,n\}$, we can always find $\xbar_1,\ldots,\xbar_n\in\{0,1\}$ so
that the partial sums $\xbar_1+\cdots +\xbar_k$ and $\xbar_{\sigma
1}+\cdots +\xbar_{\sigma k}$ differ from the unrounded values $x_1+\cdots +
x_k$ and $x_{\sigma 1}+\cdots +x_{\sigma k}$ by at most $n/(n+1)$, for $1\leq
k\leq n$. The latter bound is best possible. The proof uses an elementary
argument about flows in a certain network, and leads to a simple algorithm that
finds an optimum way to round.
\smallskip}

\bigskip
\noindent
Many combinatorial optimization problems in integers can be solved or
approximately solved by first obtaining a real-valued solution and then
 rounding
to integer values. Spencer [\cite\spenc] proved that it is always possible
to do the rounding so that partial sums in two independent orderings are
properly rounded. His proof was indirect---a~corollary of more general
results~[\cite\lsv] about discrepancies of set systems---and
 it guaranteed only that the
rounded partial sums would
differ by at most $1-2^{{-2}^n}$ from the unrounded values.
The purpose of this note is to give a more direct proof, which leads to a
sharper result.

Let $x_1,\ldots,x_n$ be real numbers and let $\sigma$ be a permutation of
$\{1,\ldots, n\}$. We will write 
$$S_k=x_1+\cdots +x_k\,,\qquad \Sigmait_k=x_{\sigma 1}+\cdots +x_{\sigma
k}\,,\qquad 0\leq k\leq n\,,$$
for the partial sums in two independent orderings. Our goal is to find integers
$\xbar_1,\ldots,\xbar_n$ such that
$$\lfloor x_k\rfloor \leq \xbar_k\leq \lceil x_k\rceil\,,$$
and such that the rounded partial sums
$$\Sbar_k=\xbar_1+\cdots +\xbar_k\,,\qquad \Sigbar_k=\xbar_{\sigma
1}+\cdots +\xbar_{\sigma k}$$
also satisfy
$$\lfloor S_k\rfloor \leq \Sbar_k\leq \lceil S_k\rceil\,,\qquad
\lfloor\Sigmait_k\rfloor\leq
\Sigbar_k\leq\lceil\Sigmait_k\rceil\,,\eqno(\ast)$$
for $0\leq k\leq n$. Such $\xbar_1,\ldots,\xbar_n$ will be called a {\it
two-way rounding\/} of $x_1,\ldots,x_n$ with respect to~$\sigma$.

\proclaim
Lemma. Two-way rounding is always possible.

\proof
We can assume without loss of generality that $S_n=m$ is an integer, by adding
an additional term and increasing~$n$ if necessary. We can also assume that
$0< x_k< 1$ for all~$k$. Construct a network with nodes
$\{s,\,a_1,\ldots,a_m,\,u_1,\ldots,u_n,\,v_1,\ldots,v_n,\,b_1,\ldots,b_m,\,t\}$
and the following arcs:\footnote*{Here and in the sequel $[a\ldt b)$ denotes
the half-open interval $\{\,x\,\vert\,a\leq x<b\,\}$. This notation, due
independently to Hoare and Ramshaw, is recommended in~[\cite\gkp].}
$$\displaylines{s\rar a_j\qquad\hbox{and}\qquad b_j\rar t\qquad\hbox{for}\qquad
1\leq j\leq m\,;\cr
\noalign{\smallskip}
u_k\rar v_k\qquad\hbox{for}\qquad 1\leq k\leq n\,;\cr
\noalign{\smallskip}
a_j\rar u_k\qquad\hbox{if}\qquad [j-1\ldt j)\cap [S_{k-1}\ldt S_k)\neq
\emptyset\,;\cr
\noalign{\smallskip}
v_{\sigma k}\rar b_j\qquad\hbox{if}\qquad
[j-1\ldt j)\cap[\Sigmait_{k-1}\ldt \Sigmait_k)\neq \emptyset\,.\cr}$$
Each arc has capacity 1. This network supports a natural flow of $m$ units, if
we send 1~unit through each arc $s\rar a_j$ and $b_j\rar t$, and $x_k$~units
through $u_k\rar v_k$; the flow in $a_j\rar u_k$ is the measure of the interval
$[j-1\ldt j)\cap[S_{k-1}\ldt S_k)$, and the flow in $v_{\sigma k}\rar b_j$ is
similar. Deleting the arcs $s\rar a_j$ defines a cut of capacity~$m$, so this
must be a minimum cut.

Since the arc capacities are integers, the max-flow/min-cut theorem implies
that this network supports an integer flow of $m$~units. Let $\xbar_k$ be the
amount that flows through $u_k\rar v_k$, for $1\leq k\leq n$, in one such flow.
Then $\xbar_k\in\{0,1\}$. If $j=\lceil S_k\rceil$ we have
$\Sbar_k=\xbar_1+\cdots +\xbar_k=$ flow into $\{u_1,\ldots,u_k\}\leq$ flow
out of $\{a_1,\ldots,a_j\}=j$, because all arcs $a_i\rar u_l$ for $l\leq k$
have $i\leq j$. 
If $j=\lfloor S_k\rfloor$ then $\Sbar_k=$ flow into $\{u_1,\ldots, u_k\}\geq$
flow out of $\{a_1,\ldots,a_j\}=j$, because all arcs $a_i\rar u_l$ for $i\leq
j$ have $l\leq k$. A~similar argument proves that $\lfloor \Sigmait_k\rfloor
\leq\Sigbar_k\leq\lceil\Sigmait_k\rceil$, hence $(\ast)$ holds.\ \pfbox

\proclaim
Corollary.
Given any fixed $k$, two-way rounding is possible with $\xbar_k=\lceil
x_k\rceil$, as well as with $\xbar_k=\lfloor x_k\rfloor$.

\proof
We may assume as before that $0<x_k<1$.
The construction in the lemma establishes a feasible flow of $x_k$~units in the
arc $u_k\rar v_k$. It is well known that the polytope of all feasible flows has
vertices whose coordinates are integers (see, for example, Application 19.2 in
Schrijver~[\cite\schrij]).
Therefore the arc $u_k\rar v_k$ is saturated in at least one
maximum flow, and it carries no flow at all in at least one other. \ \pfbox

\medskip
Incidentally, it is important to impose a capacity of 1 on the arcs $u_k\rar
v_k$ in the construction of this proof. Otherwise we might get solutions in
which $\xbar_k=2$. Condition $(\ast)$ does not by itself imply that
$\xbar_k\leq\lceil x_k\rceil$ or that $\xbar_k\geq\lfloor x_k\rfloor$.

Notice that $(\ast)$ is equivalent to the conditions
$$\vert S_k-\Sbar_k\vert < 1\qquad\hbox{and}\qquad
\vert\Sigmait_k-\Sigbar_k\vert< 1\,,\qquad \hbox{for}\quad 0\leq k\leq n\,,$$
since $\Sbar_k$ and $\Sigbar_k$ are integers. Let us say that two-way rounding
has discrepancy bounded by~$\delta$ if $\vert S_k-\Sbar_k\vert\leq\delta$ and
$\vert\Sigmait_k-\Sigbar_k\vert\leq\delta$ for all~$k$. A~slight extension of
the construction in the lemma makes it possible to prove a stronger result:

\proclaim
Theorem 1. If $S_n=m$ is an integer, the sequence $(x_1,\ldots,x_n)$ can be
two-way rounded with discrepancy bounded by $(2m+1)/(2m+2)$. 

\proof
We will prove that two-way rounding bounded by $\delta$ is possible for all
$\delta >\break
(2m+1)/(2m+2)$. Only finitely many roundings exist, so the stated
result follows by taking the limit as $\delta$~decreases to $(2m+1)/(2m+2)$.

The proof uses a network like that of the lemma, but we omit certain arcs that
would lead to discrepancies near~1. More precisely, if $\epsilon$ is any fixed
positive number $< 1/(2m+2)$, we have
$$\displaylines{a_j\rar u_k\qquad\hbox{if}\quad
[j-1+\epsilon\ldt j-\epsilon)\cap[S_{k-1}\ldt S_k)\neq\emptyset\,;\cr
\noalign{\smallskip}
v_{\sigma k}\rar b_j\qquad\hbox{if}\quad
[j-1+\epsilon\ldt j-\epsilon)\cap[\Sigmait_{k-1}\ldt
\Sigmait_k)\neq\emptyset\,.\cr}$$
We also allow these arcs to have infinite capacity. But the capacity of the
``source'' arcs $s\rar a_j$, the ``middle'' arcs $u_k\rar v_k$, and the
``sink'' arcs $b_j\rar t$ remains~1.

The minimum cut in this reduced network has size~$m$. For if any $m-1$ of the
unit-capacity arcs are cut, we will prove that we can still connect~$s$ to~$t$.
Suppose we remove $p$~source arcs, $q$~middle arcs, and $r$~sink arcs, where
$p+q+r=m-1$. We send $1-2\epsilon$ units of flow from~$s$ through each of the
$m-p$ remaining source arcs. From every~$a_j$ reached in this way, we send as
many units of flow from $a_j\rar u_k$ as the size of the interval 
$[j-1+\epsilon\ldt j-\epsilon)\cap [S_{k-1}\ldt S_k)$. Some of the flow now
gets stuck, if $u_k$ is one of the $q$~vertices for which the arc
$u_k\rar v_k$ was removed. But at most $1-2\epsilon$ units flow into
each~$u_k$, so we still have at least $(m-p-q)(1-2\epsilon)=(r+1)(1-2\epsilon)$
units of flow arriving at $\{v_1,\ldots,v_n\}$. Now consider an ``antiflow'' of
$1-2\epsilon$ units from~$t$ back through each of the $m-r$ remaining sink
arcs $b_j\rar t$. From every such~$b_j$ we send the antiflow back through
$v_{\sigma k}\rar b_j$ according to the size of $[j-1+\epsilon\ldt
j-\epsilon)\cap [\Sigmait_{k-1}\ldt \Sigmait_k)$. In this way
$(m-r)(1-2\epsilon)$ units of antiflow come from~$t$ to $\{v_1,\ldots,v_n\}$.
Each vertex~$v_k$ contains at most $x_k$~units of flow and at most $x_k$~units
of antiflow. We know that the total flow plus antiflow at $\{v_1,\ldots,v_n\}$
is at least $(r+1)(1-2\epsilon)+(m-r)(1-2\epsilon)=m+1-(2m+2)\epsilon
>m=x_1+\cdots+x_n$. Therefore some vertex~$v_k$ must contain both flow and
antiflow. And this establishes the desired link between~$s$ and~$t$.

Since $m$ is the size of a minimum cut and all capacities are integers, the
network supports an integer flow of value~$m$. Let $\xbar_k$ be the flow
from~$u_k$ to~$v_k$; we will prove that
$(\,\xbar_1,\ldots,\xbar_n)$ is a two-way rounding with discrepancy
$<\delta=1-\epsilon$. Note that
$$\vert\Sbar_k-S_k\vert< 1-\epsilon\;\Leftrightarrow\;
\lfloor S_k+\epsilon\rfloor\leq\Sbar_k\leq\lceil S_k-\epsilon\rceil\,.$$
If $j=\lceil S_k-\epsilon\rceil$ we have
$\Sbar_k=\xbar_1+\cdots+\xbar_k=$ flow into $\{u_1,\ldots,u_k\}\leq$ flow out
of $\{a_1,\ldots,a_j\}=j$, because all arcs $a_i\rar u_l$ for $l\leq k$ have
$[i-1+\epsilon\ldt i-\epsilon)\cap[S_{l-1},S_l)\neq\emptyset$, hence
$i-1+\epsilon< S_l$ and $i\leq \lceil S_l-\epsilon\rceil\leq j$. Similarly,
if $j=\lfloor S_k+\epsilon\rfloor$ we have $\Sbar_k\geq$ flow out of
$\{a_1,\ldots,a_j\}=j$, because all arcs $a_i\rar u_l$ for $i\leq j$ have
$l\leq k$. (If $l> k$ we would have $S_{l-1}\geq S_k\geq j-\epsilon\geq
i-\epsilon$, contradicting $S_{l-1}< i-\epsilon$.) A~similar proof shows 
that $\lfloor\Sigmait_k+\epsilon\rfloor
\leq\Sigbar_k\leq\lceil\Sigmait_k-\epsilon\rceil$.\ \pfbox

\smallskip
The bound of Theorem 1 is, in fact, best possible, in the sense that no better
bound can be guaranteed as a function of~$m$. 

\proclaim
Theorem 2. For all positive integers $m$ there exists a sequence of real
numbers $(x_1,\ldots,x_n)\kern-1pt$ with sum~$m$
and a permutation~$\sigma$ of $\{1,\ldots,n\}$ that
cannot be two-way rounded with discrepancy  $<(2m+1)/(2m+2)$. 

\proof
Let $n=2m+2$ and $\epsilon =1/n$. Define
$$\eqalign{%
&x_1=x_2=x_3=\epsilon\,;\qquad x_{m+3}=(2m-1)\epsilon\,;\cr
\noalign{\smallskip}
&x_{k+3}=2\epsilon\,,\;x_{k+m+3}=2m\epsilon\,,\quad\hbox{for}\quad 1\le k<
m\,;\cr
\noalign{\smallskip}
&\sigma 1=2\,,\;\sigma 2=1\,,\;\sigma 3=m+3\,,\;\sigma(2m+2)=3\,;\cr
\noalign{\smallskip}
&\sigma(2k+2)=k+3\,,\;\sigma(2k+3)=k+m+3\,,\quad\hbox{for}\quad 1\leq k<
m\,.\cr}$$
For example, when $m=4$ we have $(x_1,\ldots,x_{10})=(.1,\,.1,\,.1,\,.2,\,
.2,\,.2,\,.7,\,.8,\,.8,\,.8)$ and $(\sigma 1,\ldots,\sigma 10)=
(2,1,7,4,8,5,9,6,10,3)$. Hence
$$\eqalign{(S_1,\ldots,S_{10})&=(.1,\,.2,\,.3,\,.5,\,.7,\,.9,\,1.6,\,
2.4,\,3.2,\,4.0)\,,\cr
\noalign{\smallskip}
(\Sigmait_1,\ldots,\Sigmait_{10})&=(.1,\,.2,\,.9,\,1.1,\,1.9,\,2.1,\,2.9,\,
3.1,\,3.9,\,4.0)\,.\cr}$$
We will prove that this sequence and permutation cannot be two-way rounded with
discrepancy less than $(2m+1)/(2m+2)=0.9$; the same proof technique will work
for any $m\geq 1$.

The main point is that whenever $S_k$ or $\Sigmait_k$ has the form $l\pm 0.1$
where $l$ is an integer, it must be rounded to~$l$ in order to keep the
discrepancy small. This forces $\Sbar_1=\Sigbar_1=0$,
$\Sigbar_3=\Sigbar_4=1$, $\Sigbar_5=\Sigbar_6=2$, $\Sigbar_7=\Sigbar_8=3$,
$\Sigbar_9=4$, hence $\xbar_1=\xbar_2=\xbar_3=\xbar_4=
\xbar_5=\xbar_6=0$. But then
$\Sbar_6=\xbar_1+\cdots +\xbar_6=0$ differs by 0.9 from~$S_6$.\ \pfbox

\smallskip
Although Theorem 2 proves that Theorem 1 is ``optimal,'' we can do still better
if $m$ is greater than~${1\over 2}n$, because we can replace each~$x_k$ by
$1-x_k$. This replaces $m$ by $n-m$, and the bound on discrepancy decreases to
$(2n-2m+1)/(2n-2m+2)$. Then we can restore the original $x_k$ and change
$\xbar_k$ to $1-\xbar_k$. This computation preserves $\vert S_k-\Sbar_k\vert$
and $\vert\Sigmait_k-\Sigbar_k\vert$, so it preserves the discrepancy.

Further improvement is also possible when $m=\lfloor n/2\rfloor$, if we look at
the construction closely. The following theorem gives a uniform bound in terms
of~$n$, without any assumption about the value of $x_1+\cdots+x_n$.

\proclaim
Theorem 3. Any sequence $(x_1,\ldots,x_n)$ and permutation $(\sigma 1,\ldots,
\sigma n)$ can be two-way rounded with discrepancy bounded by $n/(n+1)$.

\proof
We will show in fact that the discrepancy can always be bounded by $(n-1)/n$,
when $x_1+\cdots +x_n=m$ is an integer. The general case follows from this
special case if we set $x_{n+1}=\lceil S_n\rceil -S_n$ and increase~$n$ by~1.

If $2m+2\leq n$ or $2n-2m+2\leq n$, the result follows from Theorem~1 and
possible complementation. Therefore we need only show that a discrepancy of at
most $(n-1)/n$ is achievable when $m=\lfloor n/2\rfloor$.

Consider first the case $n=2m+1$. We use the network in the proof of Theorem~1,
but now we allow $\epsilon$ to be any  number $< 1/n$. Suppose, as in
the former proof, that we can disconnect~$s$ from~$t$ by deleting $p$~source
arcs, $q$~middle arcs, and $r$~sink arcs,
where $p+q+r=m-1$. Let $q$ be minimum over all such ways
to disconnect the network. We construct flows and antiflows as before, and we
say that $x_k$ is {\it green\/} if $v_k$~contains positive flow, {\it red\/} if
$v_k$~contains positive antiflow. No $x_k$ is both green and red, since there
is no path from~$s$ to~$t$. The previous proof showed that there are at least
$(r+1)(1-2\epsilon)$ units of green flow and $(m-r)(1-2\epsilon)$ units of red
flow, hence there are at least $m+1-(2m+2)\epsilon$ units of flow altogether.
If we can raise this lower bound by~$\epsilon$, we will have a contradiction,
because $m+1-(2m+1)\epsilon> m$.

Suppose $q> 0$, and let $u_k\rar v_k$ be a middle arc that was deleted. At
most two arcs emanate from~$v_k$ in the network. Since $q$ is minimum, there
must in fact be two; otherwise we could restore $u_k\rar v_k$ and delete a
non-middle arc. The two arcs from~$v_k$ must be consecutive, from $v_k\rar b_j$
and $v_k\rar b_{j+1}$, say. Furthermore the arcs $b_j\rar t$ and $b_{j+1}\rar
t$ have not been cut. If $k=\sigma l$ we have $\Sigmait_{l-1}< j-\epsilon$
and $\Sigmait_l> j+\epsilon$. Our lower bound on antiflow can now be raised
by~$2\epsilon$, because it was based on the weak assumption that no antiflow
runs back from $[j-\epsilon\ldt j+\epsilon)$. This improved lower bound leads
to a contradiction; hence $q=0$.

Divide the interval $[0\ldt m)$ into $3m$~regions, namely ``tiny left'' regions
of the form $[j-1\ldt j-1+\epsilon)$, ``inner'' regions of the form
$[j-1+\epsilon\ldt j-\epsilon)$, and ``tiny right'' regions of the form
$[j-\epsilon\ldt j)$, for $1\leq j\leq m$. If we color the points of
$[S_{k-1}\ldt S_k)$ with the color of~$x_k$, our lower bound
$(r+1)(1-2\epsilon)$ for green flow was essentially obtained by noting that
$m-p=r+1$ of the inner regions are purely green. Similarly, if we color the
ponts of $[\Sigmait_{k-1}\ldt \Sigmait_k)$ with the color of~$x_{\sigma k}$,
our lower bound for red flow was obtained by noting that $m-r=p+1$ inner
regions in this second coloring are purely red. Notice that there is complete
symmetry between red and green, because we can invert the network and
replace~$\sigma$ by~$\sigma^{-1}$.

Call an element $x_k$ {\it large\/} if it exceeds $1-\epsilon$. If any $x_k$ is
large, the interval $[S_{k-1}\ldt S_k)$ occupies more than $\epsilon$~units
outside of an inner region; this allows us to raise the lower bound
by~$\epsilon$ and obtain a contradiction. Therefore no element is large. It
follows that no element~$x_k$ can intersect more than 2~tiny regions, 
when $x_k$~is
placed in correspondence with $[S_{k-1}\ldt S_k)$ or with 
 $[\Sigmait_{k-1}\ldt \Sigmait_k)$.

Let's look now at the $2m$ tiny regions. Each of them must contain at least
some red in the first coloring; otherwise we would have at least
$(p+1)(1-2\epsilon)$ red units packed into at most $2m-1$ tiny regions and
$p$~inner regions, hence $(p+1)(1-2\epsilon)\leq (2m-1)\epsilon
+p(1-2\epsilon)$, contradicting $\epsilon< 1/n$. This means there must be at
least $m+1$ red elements~$x_k$, since no red element is large and since
$m$~non-large red intervals can intersect all the tiny regions only if they
also cover all the inner regions (at least one of which is green).
Similarly, there must be at least $m+1$ green
elements. But this is impossible, since there are only $2m+1$ elements
altogether. Therefore the network has minimum cut size~$m$, and the rest of the
proof of Theorem~1 goes through as before.

Now suppose $n=2m$. Then we can carry out a similar argument, but we need to
raise the lower bound by~$2\epsilon$. Again we can assume that $q=0$. We can
also show without difficulty that there cannot be {\it two\/} large elements.
When $n=2m$ the argument given above shows that at least $2m-1$ of the tiny
regions must contain some red, in the first coloring.

Suppose there are only $m-1$ red elements. Then, in the first coloring, $m-2$
of them intersect 2~tiny intervals and the other is large and intersects~3; we
have raised the red lower bound by~$\epsilon$. 
But $(p+1)(1-2\epsilon)+\epsilon$ red units cannot be packed into $2m-1$ tiny
regions and $p$~inner regions, because $(p+1)(1-2\epsilon)+\epsilon >
(n-1)\epsilon +p(1-2\epsilon)$. 

A symmetrical argument shows that there cannot be only $m-1$ green elements.
Therefore exactly $m$ elements are red and exactly $m$ are green. Suppose no
element is large. Then we have at least one purely green tiny interval in the
first coloring and at least one purely red tiny interval in the
second---another contradiction. Thus, we may assume that there is one large red
element, and that the $2m$~tiny intervals in the first coloring 
contain a total of less than $\epsilon$~units of green. In particular, each of
them contains some red. Either the first interval $[0\ldt \epsilon)$ or the
last interval $[m-\epsilon \ldt m)$ is intersected by a non-large red element,
which intersects at most $\epsilon$~units of space in tiny intervals. The other
$m-1$ red elements intersect at most $2\epsilon$~units of tiny space each, so
at most $(2m-1)\epsilon$ such units are red. This final contradiction completes
the proof. \ \pfbox

\smallskip
The result of Theorem 3 is best possible, because we can easily prove (as in
Theorem~2) that the values
$$x_1={1\over n+1}\,,\qquad x_k=\cases{(n-1)/(n+1)\,,&$k$ even, $2\leq k\leq
n$\cr
\noalign{\smallskip}
2/(n+1)\,,&$k$ odd, $3\leq k\leq n$\cr}$$
and a ``shuffle'' permutation that begins 
$$\sigma k=\left\{\,\vcenter{\halign{$#$\hfil\quad&#\hfil\quad&#\hfil\cr
2k-1&for $1\leq 2k-1\leq n$,&$n$ odd\cr
\noalign{\smallskip}
2k&for $1\leq 2k\leq n$,&$n$ even\cr}}\right.$$
cannot be two-way rounded with discrepancy less than $n/(n+1)$.

So far we have discussed only worst-case bounds. But a particular two-way
rounding problem, defined by values $(x_1,\ldots,x_n)$ and a permutation
$(\sigma 1,\ldots,\sigma n)$, will usually be solvable with smaller discrepancy
than guaranteed by Theorems~1 and~3. 
A~closer look at the construction of Theorem~1
leads to an efficient algorithm that finds the best possible discrepancy in any
given case.

\proclaim Theorem 4.
Let $\epsilon$ be any positive number.
There exists a solution with discrepancy less than $1-\epsilon$ to a given
two-way rounding problem if and only if the network constructed in the proof of
Theorem~1 supports an integer flow of value~$m$.

\proof
The final paragraph in the proof of Theorem 1 demonstrates the ``if'' half.
Conversely, suppose $\xbar_1,\ldots,\xbar_n$ is a
solution with discrepancy $<1-\epsilon$. If $\xbar_k=1$, let $j=\bar{S}_k$.
Then $j-1=\bar{S}_{k-1}$, so the condition
$\vert\bar{S}_{k-1}-S_{k-1}\vert <1-\epsilon$
implies $S_{k-1}<j-\epsilon$. Also $\vert\bar{S}_k-S_k\vert <1-\epsilon$
implies $S_k>j-1+\epsilon$. Therefore there is an arc $a_j\rightarrow u_k$.
Similarly, there is an arc $v_{\sigma k}\rightarrow b_j$ when $\xbar_{\sigma
k}=1$ and 
$j=\overline{\Sigmait}_k$. So the network supports an integer
flow of value~$m$. \ \pfbox

\smallskip
In other words, the optimum discrepancy $\delta =1-\epsilon$ is obtained when
$\epsilon$ is just large enough to reduce the network to the point where no
$m$-unit flow can be sustained,
if $\delta \ge {1\over 2}$. We can in fact find an optimum rounding as follows:
Let 
$$f(j,k)=\min(j-S_{k-1},S_k-j+1)$$
be the ``desirability'' of the arc $a_j\rar u_k$, and
$$g(j,\sigma k)=\min(j-\Sigmait_{k-1},\Sigmait_k-j+1)$$
the desirability of $v_{\sigma k}\rar b_j$. (Thus the arcs $a_j\rar u_k$,
$v_{\sigma k}\rar b_j$ are included in the network of Theorem~1 if and only if
their desirability is greater than~$\epsilon$.) Sort these arcs by
desirability, and add them one by one to the initial arcs $\{s\rar a_j,u_k\rar
v_k,b_j\rar t\}$ until an integer flow of $m$~units is possible. Then let
$\xbar_k$ be the flow in $u_k\rar v_k$, for all~$k$; this flow has discrepancy
equal to~1 minus the desirability of the last arc added, and no smaller
discrepancy is possible.

Notice that the arc $a_j\rar u_k$ has desirability $>{1\over 2}$ if and only if
$S_{k-1}<j-{1\over 2}<S_k$, so at most $m$ such arcs are present.
 If all $x_k$ lie between~0 and~1, at most $m+n-1$ arcs 
of the form $a_j\rar u_k$ will have positive
desirability, since both $a_{j-1}\rar u_k$ and $a_j\rar u_k$ will be desirable
iff $S_{k-1}<j<S_k$. 

The following simple algorithm turns out to be quite efficient, assuming that
$m\le {1\over 2}n$: Begin with the network consisting of arcs $\{s\rar a_j,
u_k\rar v_k, b_j\rar t\}$ for $1\le j\le m$ and $1\le k\le n$, plus any
additional arcs of desirability $>{1\over 2}$. Call an arc $a_j\rar u_k$ or
$v_{\sigma k}\rar b_j$ ``special'' if its desirability lies between
$1/\min(2m+2,n)$ and~${1\over 2}$, inclusive; fewer than $2m+2n$ arcs are
special. 
Then, for $j=1,\ldots,m$, send one unit of flow from~$a_j$ to~$t$ along an
``augmenting path,'' using the well-known algorithm of Ford and Fulkerson
[\cite\forful,
pp.\ 17--19] but specialized for unit-capacity arcs. In other words, construct
a breadth-first search tree from~$a_j$ until encountering~$t$; then choose a
path from~$a_j$ to~$t$ and reverse the orientation of all arcs on that path.
If $t$ is not reachable from~$a_j$, add special arcs to the network, in order
of decreasing desirability, until $t$ is reachable.

The running time of this algorithm is bounded by $O(mn)$ steps, but in practice
it runs much faster on random data. For example, Tables~1 and~2 show the 
results
of various tests when the input permutation~$\sigma$ is random and when
the values
$(x_1,\ldots,x_n)$ are selected as follows: Let $y_1,\ldots,y_n$ be
independent uniform integers in the range $1\le y_k\le N$, where $N$ is a large
integer (chosen so that arithmetic computations will not exceed 31~bits).
Increase one or more
 of the~$y$'s by~1, if necessary, until $y_1+\cdots +y_n$ is a multiple
of~$m$; then set $x_k=y_k/d$, where $d=(y_1+\cdots +y_n)/m$. Reject
$(x_1,\ldots,x_n)$ and start over, if some $x_k\ge 1$. (In practice, rejection
occurs about half the time when $m={1\over 2} n$, but almost never when $m\ll
{1\over 2} n$.)

Table 1 shows the optimum discrepancies found, and Table~2 shows the running
time in memory references or ``mems'' [\cite\sgb,
 pp.\ 464--465] divided by~$n$. All
entries in these tables are given in the form $\mu\pm\sigma$, where $\mu$ is
the sample mean and $\sigma$ is an estimate of the standard deviation; more
precisely, $\sigma$~is the square root of an unbiased estimate of the variance.
The number of test runs $t(n)$ for each experiment was $10^6\!/n$; thus,
$10^5$~runs were made for each~$m$ when $n=10$, but only 10~runs were made for
each~$m$ when $n=10^5$. The actual confidence interval for the tabulated
$\mu$~values is therefore approximately
$2\sigma/\sqrt{t(n)}=.002\sigma\sqrt{n}$. 

\vfill\eject
\centerline{Table 1.\quad Empirical optimum discrepancies}
$$\vcenter{\halign{$#$\hfil\quad
&$#$\hfil\quad&$#$\hfil\quad&$#$\hfil\quad&$#$\hfil\quad&$#$\hfil\cr
&\hfil m=1&\hfil m=2&\hfil m=\lfloor\lg n\rfloor&\hfil m=\lfloor\sqrt{n}\rfloor
&\hfil m={1\over 2}n\cr
\noalign{\smallskip}
n=10&.566\pm .06&.619\pm .07&.627\pm .07&.627\pm .07&.622\pm .08\cr
n=100&.537\pm .02&.575\pm .03&.664\pm .03&.710\pm .03&.759\pm .02\cr
n=1000&.513\pm .007&.527\pm.01&.582\pm.01&.662\pm.02&.794\pm.02\cr
n=10000&.504\pm.002&.509\pm.003&.535\pm.005&.612\pm.01&.818\pm.01\cr
n=100000&.502\pm.001&.503\pm.001&.513\pm.002&.570\pm.005&.838\pm.007\cr
}}$$

\bigskip
\centerline{Table 2.\quad Empirical running time, in mems/{\it n}}
$$\vcenter{\halign{$#$\hfil\quad
&$#$\hfil\quad&$#$\hfil\quad&$#$\hfil\quad&$#$\hfil\quad&$#$\hfil\cr
&\hfil m=1&\hfil m=2&\hfil m=\lfloor\lg n\rfloor&m=\lfloor\sqrt{n}\rfloor
&\hfil m={1\over 2}n\cr
\noalign{\smallskip}
n=10&\phantom{.}10\pm 4&\phantom{.}19\pm 6&\phantom{.}27\pm 8&\quad 27\pm 8%
&\phantom{1}37\pm 11\cr
n=100&2.9\pm 1.3&\phantom{1.}6\pm 2&\phantom{.}18\pm 5%
&\quad 29\pm 7&\phantom{1}76\pm 15\cr
n=1000&0.9\pm 0.5&1.9\pm 0.7&8.5\pm 2.2&\quad 25\pm 6&152\pm 32\cr
n=10000&0.3\pm 0.2&0.6\pm 0.2&3.6\pm 0.8&\quad 22\pm 7&289\pm 49\cr
n=100000&0.1\pm 0.1&0.2\pm 0.1&1.4\pm 0.4&\quad 17\pm 4&540\pm 72\cr
}}$$

Notice that when $m\ll n$, the optimum discrepancy is nearly~${1\over 2}$.
Indeed, this is obvious on intuitive grounds: When $n$ is large, approximately
$\epsilon n$~values of~$k$ will have $S_k$ within ${1\over 2}\epsilon$ of
$\{{1\over 2}\,,\,{3\over 2}\,,\ldots,m-{1\over 2}\}$, and approximately
$\epsilon^2n$ will also have equally good values~$\Sigmait_{\sigma^{-1}k}$.
So we are essentially looking for a perfect matching in a bipartite graph with
$m$~vertices in each part and $\epsilon^2n$ edges. For fixed $m$ as $n\rar
\infty$, the matching will exist when $\epsilon^2n$ is sufficiently large,
hence the mean optimum discrepancy is ${1\over 2}+O(n^{-{1\over 2}})$.

However, the behavior of the mean optimum discrepancy when $m={1\over 2}n$ is
not clear. It appears to approach~1, but quite slowly, perhaps as $1-c/\log n$.

When $n$ is fixed and $m$ varies, the mean optimum discrepancy is not maximized
when $m={1\over 2}n$.  For example, when $n=10$, Table~1 shows that it is .622
when $m=5$ but .627 when $m=3$. 

The running times shown in Table~2 do not include the work of constructing the
network or sorting the special arcs by desirability. Those operations are
easily analyzed, and in practice they take $am+bn$ steps for some constants~$a$
and~$b$, because a straightforward bucket sort is satisfactory for this
application. Therefore only the running time of the subsequent flow
calculations is of interest.

The average running time to compute the flows appears to be $o(n)$ when
$m\le\sqrt{n}$, and approximately proportional to~$n^{1.3}$ when $m={1\over
2}n$. So it is much less than the obvious upper bound~$mn$ of the
Ford-Fulkerson scheme. The author tried to obtain still faster results by using
more sophisticated max-flow algorithms, but these ``improved'' algorithms
actually turned out to run more than an order of magnitude slower.

For example, the algorithm of Dinits, as improved by Karzanov and others, seems
at first to be especially well suited to this application because the network
of Theorem~1 is ``simple'' in the sense discussed by Papadimitriou and
Steiglitz [\cite\papst, pp.\ 212--214]: Every internal vertex 
has in-degree~1 or out-degree~1, hence edge-disjoint paths are vertex-disjoint
and the running time with unit-capacity arcs is $O(\,\vert V\vert^{1/2}\,\vert
A\vert\,)=O(n^{3/2})$. Using binary search to find the optimum number of
special arcs gives us a guaranteed worst-case performance of $O\bigl(\min
(m,n^{1/2})n\log n\bigr)$. Unfortunately, in practice the performance of that
algorithm actually matches this worst-case estimate, even on random data. For
example, when $m={1\over 2}n$ the observed running time in mems/$n$ was
$15284\pm 2455$ when $n=10^4$, compared to $289\pm 49$ by the simple algorithm.
Each flow calculation consumed more than $1000n$ mems, and binary search
required $\lceil\lg 2n\rceil =14$ flow calculations to be carried out.

When modern preflow push/relabel algorithms are specialized to unit-capacity
networks of the type considered here, they behave essentially like the Dinits
algorithm and are no easier to implement (see Goldberg, Plotkin, and Vaidya
[\cite\gpv]). Such algorithms do allow networks to change dynamically by adding
arcs from~$s$ and/or deleting arcs to~$t$ (see Gallo, Grigoriadis, and Tarjan
[\cite\ggt]);
but our application requires adding or deleting special arcs in the {\it
middle\/} of the network, so the 
techniques of~[\cite\ggt] do not apply. Thus the simple
Ford-Fulkerson algorithm 
seems to be a clear winner for this application, in spite of a
lack of performance guarantees.

How complex can the networks of Theorem 1 be? If we have any bipartite graph
with $m$~vertices in each part and with $n$~edges, and if every edge can be
extended to a perfect matching, then we can find real numbers
$(x_1,\ldots,x_n)$ in the range $0<x_k\le 1$ and a permutation $(\sigma
1,\ldots, \sigma n)$ such that $x_1+\cdots +x_n=m$ and the two-way roundings
are in one-to-one correspondence with the perfect matchings of the given graph.
For we can take $(x_1,\ldots,x_n)=t_1\alpha_1+\cdots +t_n\alpha_n$ where
$t_1+\cdots +t_n=1$ and $\alpha_k$ is the characteristic vector of a perfect
matching that uses edge~$k$. The sum of $x_k$ 
over all the edges touching any vertex
is~1. Represent an edge from~$u$ to~$v$ by the ordered pair $(u,v)$, and label
the edges $1,\ldots,n$ in lexicographic order of these pairs; then define the
permutation $\sigma 1,\ldots,\sigma n$ by lexicographic order of the dual pairs
$(v,u)$. It follows that if $k$ is the final edge for vertex~$j$ 
in the first part, we have $S_k=j$; and if $\sigma k$ is the final edge for
vertex~$j$ 
in the second
part, we have $\Sigmait_k=j$. The correspondence between matchings and
roundings is now evident.

This construction shows that the networks arising in Theorem~1 are general
enough to mimic the networks that arise in bipartite matching problems, but
only when the bipartite graphs contain no unmatchable edges; and the corollary
preceding Theorem~1 shows that the latter restriction cannot be removed. This
restriction on network complexity
might account for the excellent performance we obtain with the
simple Ford-Fulkerson algorithm.

If the capacity constraint on $u_k\rar v_k$ is removed, our network becomes
equivalent to a network for bipartite matching, in which we want to match
$\{a_1,\ldots,a_m\}$ to $\{b_1,\ldots,b_m\}$ through edges $a_j\rbar b_{j'}$
whenever $a_j\rar u_k$ and $v_k\rar b_{j'}$. The problem of finding the {\it
best\/} such match, when the edge $a_j\rbar b_{j'}$ is ranked by the minimum of
the desirabilities $f(j,k)$ and $g(j',k)$, is then a {\it bottleneck assignment
problem\/} [\cite\edful, \cite\forful].
(Open question: Is there a nice way to characterize all bottleneck 
assignment problems that arise from two-way rounding problems in this
manner?)

The problem of optimum two-way rounding is, however, more general than the
bottleneck assignment problem, because the unit capacity constraint on $u_k\rar
v_k$ is significant. Consider, for example, the case $n=7$, $m=3$,
$(x_1,\ldots,x_7)={1\over 28}(8,8,24,11,11,11,11)$, $(\sigma 1,\ldots,\sigma 7)
=(2,1,3,5,4,7,6)$. Then $(S_1,\ldots,S_7)=(\Sigmait_1,\ldots,\Sigmait_7)=
{1\over 28}(8,16,40,51,62,\allowbreak
73,84)$, and the arcs $\{a_j\rar u_k,v_k\rar b_j\}$
ranked by desirability are
$$\vcenter{\halign{$#$\hfil\qquad&\hfil#\cr
a_3\rar u_{6},\ v_7\rar b_3&desirability $=\min\!\left(\,{22\over 28}\,,\,
{17\over 28}\,\right)={17\over 28}$\cr
\noalign{\smallskip}
a_1\rar u_2,\ v_1\rar b_1,\ a_2\rar u_4,\ v_5\rar b_2&desirability ${16\over
28}$\cr
\noalign{\smallskip}
a_1\rar u_3,\ v_3\rar b_1,\ a_2\rar u_3,\ v_3\rar b_2&desirability ${12\over
28}$\cr
\noalign{\smallskip}
a_3\rar u_7,\ v_6\rar b_3&desirability ${11\over 28}$\cr
\noalign{\smallskip}
a_1\rar u_1,\ v_2\rar b_1&desirability ${8\over 28}$\cr
\noalign{\smallskip}
a_3\rar u_5,\ v_4\rar b_3&desirability ${6\over 28}$\cr
\noalign{\smallskip}
a_2\rar u_5,\ v_4\rar b_2&desirability ${5\over 28}$\cr}}$$
Thus the edges $a_j\rbar b_{j'}$ ranked by desirability are
$$\vcenter{\halign{$#$\hfil\qquad&$#$\hfil\cr
a_1\rbar b_1,\ a_1\rbar b_2,\ a_2\rbar b_1,\ a_2\rbar b_2%
&\left(\,{12\over 28}\ {\rm via}\ 
u_3,v_3\right)\cr
\noalign{\smallskip}
a_3\rbar b_3
&\left(\,{11\over 28}\ {\rm via}\ u_6,v_6 \ {\rm or} \ u_7,v_7\right)\cr
\noalign{\smallskip}
a_1\rbar b_1
&\left(\,{8\over 28}\ {\rm via}\ u_1,v_1 \ {\rm or} \ u_2,v_2\right)\cr
\noalign{\smallskip}
a_2\rbar b_3,\ a_3\rbar b_2
&\left(\,{6\over 28}\ {\rm via}\ u_4,v_4 \ {\rm or} 
\ u_5,v_5\right)\cr
\noalign{\smallskip}
a_2\rbar b_2
&\left(\,{5\over 28}\ {\rm via}\ u_4,v_4 \ {\rm or} \ u_5,v_5\right)\cr
}}$$
The bottleneck assignment problem is solved by matching  $a_1
\rbar b_1$, $a_2\rbar b_2$, and $a_3\rbar b_3$ with desirability 
$\min\left(\,{12\over 28}\,,\,{12\over 28}\,,\,{11\over 28}\,\right)={11\over
28}$. But this matching does not correspond to a valid two-way rounding because
it uses the intermediate arc $u_3\rar v_3$ twice; it rounds $x_3$ to~2 and
$x_6$ (or~$x_7$) to~1. The optimum two-way rounding uses another route
from~$a_1$ to~$b_1$ and has desirability $\min\left(\,{8\over 28}\,,\,{12\over
28}\,,\,{11\over 28}\,\right)={8\over 28}$, discrepancy $1-{8\over 28}={20\over
28}$; it rounds $x_1$ (or~$x_2$), $x_3$, and~$x_6$ (or~$x_7$) to~1, the other
$x$'s to~0.

In closing, we note that a conjecture of J\'ozsef Beck [\cite\lsv, \cite\spenc]
remains a
fascinating open problem: Is there a constant~$K$ such that three-way rounding
is always possible with discrepancy at most~$K$? $\bigl($In 
three-way rounding the
partial sums are supposed to be well approximated with respect to a third
permutation $(\tau 1,\ldots,\tau n)$, in addition to $(1,\ldots,n)$ and
$(\sigma 1,\ldots, \sigma n)$.$\bigr)$\ It suffices [\cite\lsv, \cite\spenc]
to prove this when
$x_k={1\over 2}$ for all~$k$.

Can any of the methods of this paper be extended to find better bounds on the
discrepancy of arbitrary set systems 
(or at least of set systems more general than those for two-way rounding), in
the sense of~[\cite\spenc]?

\bigskip\noindent
{\bf Acknowledgments.}
I wish to thank Joel Spencer for proposing the problem and for showing me a
simple construction that forces discrepancy $n/(n+1)$. Thanks also to Noga
Alon, Svante Janson, and Serge Plotkin
 for several stimulating discussions as I was working out the
solution described above. Shortly after I~had proved Theorems 1--3, a~somewhat
similar construction was found independently by Jacek Ossowski, who described
it in terms of common systems of distinct representatives instead of network
flows; see {\S}9.2 in~[\cite\mirsk].

\vfill\eject

{\bf References}

\medskip
\bib
\phantom1[\cite\edful]\enspace
Jack Edmonds and D. R. Fulkerson, ``Bottleneck extrema,'' {\sl Journal of
Combinatorial Theory\/ \bf 8} (1970), 299--306.

\medskip
\bib
\phantom1[\cite\forful]\enspace
L. R. Ford, Jr., and D. R. Fulkerson, {\sl Flows in Networks\/} (Princeton
University Press, 1962).

\medskip
\bib
\phantom1[\cite\ggt]\enspace
Giorgio Gallo, Michael D. Grigoriadis, and Robert E. Tarjan, ``A~fast
parametric maximum flow algorithm and applications,'' {\sl SIAM Journal on
Computing\/ \bf 18} (1989), 30--55.

\medskip
\bib
\phantom1[\cite\gpv]\enspace
Andrew V. Goldberg, Serge A. Plotkin and Pravin M. Vaidya, ``Sublinear-time
parallel algorithms for matching and related problems,'' {\sl Journal of
Algorithms\/ \bf 14} (1993), 180--213.

\medskip
\bib
\phantom1[\cite\gkp]\enspace
Ronald L. Graham, Donald E. Knuth, and Oren Patashnik, {\sl Concrete
Mathematics\/} (Addison\kern.05em--Wesley, 1989).

\medskip
\bib
\phantom1[\cite\sgb]\enspace
Donald E. Knuth, {\sl The Stanford GraphBase\/} (ACM Press, 1994).

\medskip
\bib
\phantom1[\cite\lsv]\enspace
L. Lov\'asz, J. Spencer and K. Vesztergombi, ``Discrepancy of set-systems
and matrices,'' {\sl European Journal of Combinatorics\/ \bf 7} (1986),
151--160. 

\medskip
\bib
\phantom1[\cite\mirsk]\enspace
L. Mirsky, {\sl Transversal Theory\/} (Academic Press, 1971).

\medskip
\bib
\phantom1[\cite\papst]\enspace
Christos H. Papadimitriou and Kenneth Steiglitz, {\sl Combinatorial
Optimization\/}
 (Pren\-tice-Hall, 1982).

\medskip
\bib
[\cite\schrij]\enspace
Alexander Schrijver, {\sl Theory of Linear and Integer Programming\/}
(Wiley, 1986).

\medskip
\bib
[\cite\spenc]\enspace
Joel Spencer, {\sl Ten Lectures on the Probabilistic Method}, CBMS-NSF Regional
Conference Series in Applied Mathematics, number 52 (Philadelphia: SIAM, 1987),
Lecture~5.

\bye